\documentclass{article}
\usepackage{amsmath}
\usepackage[english]{babel}
\usepackage[utf8]{inputenc}
\usepackage{graphicx}
\usepackage{soul}
\usepackage{epstopdf}
\usepackage[bottom=2.5cm,left=3cm,right=3cm,top=3cm]{geometry}

\begin{document}

\begin{center}
\section*{\LARGE Discrete  two-contour system with co-directional movement}

\textit{Yashina M.V.$^{1,2}$, Tatashev A.G.$^{1,2}$, Fomina M.J.$^{1}$}\\
\textit{$^{1}~$Moscow Automobile and Road Construction State Technical University (MADI)}\\
\textit{$^{2}~$Moscow Technical University of Communications and Informatics (MTUCI)}\\
e-mail:mv.yashina@madi.ru
\end{center}

{\bf Abstract.} {\normalsize The paper studies a discrete dynamical system, which belongs to the class of contour systems developed by A.P~Buslaev. The system contains two closed contours. There are $n$ cells and a group of particles at each contour. This group is called a cluster. The particles of a cluster are in adjacent cells and move simultaneously. At each discrete moment any particle moves onto a cell in the direction of movement. There are two common points of contours called nodes.
Delays of clusters occur at the nodes. These delays are due to that two particles cannot cross the same node simultaneously. Earlier this system were studied under the assumptions that the number of particles does not depend on the cluster. This paper studies a system with different lengths of clusters. Analytical results are obtained for
limit cycles of the system and the spectrum of average cluster velocities.

\section*{1. Introduction}

\hskip 15pt  The paper studies a discrete dynamical system, which belongs to the class of contour systems developed by A.P~Buslaev. This class was introduced to develop network transport models such that analytical results can be developed for these models.

In [1]--[7], analytical results have been obtained for mathematical traffic models such that, in these models, particles moves on a one-dimensional lattice. These models can be interpreted in terms of cellular automata, [8], or exclusive processes, [9].

The concept of traffic model was introduced in [10]. In a discrete model,
a cluster is a group of particles located in adjacent cells and moving simultaneously.
In a continuos model, a cluster is a moving segment.

The concept of contour networks (Buslaev networks) was introduced in [11]. A contour network contains closed contours, in which particles or clusters move under prescribed rules. Adjacent contours contain common points called nodes.  Delays of clusters occur at nodes. The delays are due to that more than one cluster cannot to pass through the node simultaneously. In [11]--[24] analytical results were obtained for contour network.

In [15]--[17] a contour network with two contours and two nodes were studied. The nodes divide each contour into two parts of different lengths. The length of clusters
is the same.  Two versions of the system were considered. In one of these versions
one cluster moves counter-clockwise and the other cluster moves clockwise (co-directional movement). In the other version both the clusters move clockwise (one-directional movement).

This paper studies a generalization of the system, considered in [15], under the
assumption that the movement is co-directional.

\section*{2. Description of system}

\hskip 15pt Let the system contain two closed contours,  Fig.~1. They are the contour~1 and the contour~2. There are $n$ cells in any contour. The indices of the cells are $0,1,2,\dots,n-1.$ At each discrete moment $t=0,1,2,\dots,$ each cluster moves onto one cell in the direction of movement if
the cells are numbered in the direction of movement modulo~$n.$ There is a cluster in contour~$i$ (cluster $i),$
and there are $l_i<n$ particles in the cluster $i,$ $i=1,2.$ At any moment the particles of each cluster are located in adjacent cells. There are two common points of contours. These points are called nodes. The node~1 is located between the cells $n-1$ and 0 for each contour. The node~2 is located between the cells $d-1$ and $d,$ $d\le n/2$ for each contour. We say that a cluster is {\it at the node 1 (at the node~2)} if the leading cell of the cluster is located in the cell $n-1$ (in the cell $d-1).$ We say that a cluster ({\it occupies a node}) if the node is between two particles of this cluster.
Delays occur if two particles try to cross the node simultaneously. If, at time $t,$ a cluster $i$ is at a node, and the cluster $j\ne i$ occupies this node, then the cluster $i$ does not move and, at time $t+1,$ its particles are located in the cells in that the particles are located at time $t.$ If two clusters are at the same node simultaneously,
then only the cluster~1 moves.

\begin{figure}[ht!]
\centerline{\includegraphics[width=200pt]{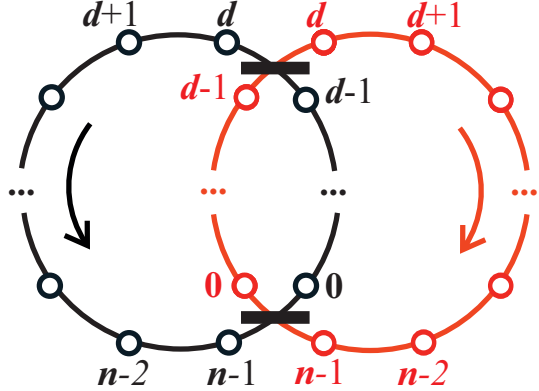}}
\caption{A two-contour system}
\end{figure}

The state of the system at time $t$ is a vector $x(t)=(x_1(t),x_2(t)),$ where $x_i(t)$ is the index of cell occupied by the leading particle of  the cluster $i,$ $i=1,2.$

The state of the system is {\it admissible} if no node is occupied by both the clusters. The {\it initial state} $x(0)$ is prescribed. This state must be admissible.

 Without loss of generality we assume that $l_1\le l_2.$

\section*{3. Limit cycles. Average velocity of cluster. States of free movement and collapse}

 \hskip 15pt The state space of the system is finite.  From any moment, system states are repeated ({\it limit cycles}).

Suppose that at time interval $(0,t)$ the cluster $i$ passes the distance
$H(t).$ The limit
\begin{equation}\label{eq1}
v_i=\lim\limits_{t\to \infty}\frac{H_i(t)}{t},
\end{equation}
is called the average velocity of the cluster $i,$ $i=1,2.$ The limit (1) exists and is equal to $A_i/T,$ where $A_i$ is the distance such that the cluster~$i$ passes this distance on the realized limit cycle, and $T$ is the period of cycle.

If, at any time $t\ge t_0$ any particle moves, then, from the moment $t_0,$  the system {\it is in a state of free movement.}

If at any moment $t\ge t_0$ no particle moves, then, from the moment  $t_0,$ the system {\it is in the state of collapse.} If the system results in the state of collapse, then $v_1=v_2=0.$

\section*{4. Some lemmas}

\hskip 18pt In this section, we shall prove some lemmas.
 \vskip 3pt
{\bf Lemma 1.} {\it If the inequality
$$l_1+l_2>n,\eqno(2)$$
hplds, then the system cannot be in a state of free movement.
\vskip 3pt
 Proof.} If, on a limit cycle, the system is in a state of free movement, then the period is equal to~$n,$ and, on the cycle, $l_1+l_2$ particles cross each node. However it is impossible if (2) holds. Indeed, two particles cannot cross the same node simultaneously. Lemma~1 has been proved.
\vskip 3pt
{\bf Lemma 2.} {\it If, on a limit cycle, the system is not in a state of free movement nor collapse, then the system results on this cycle in at least one of the states $(l_1,0),$ $(0,l_2),$ $(d+l_1,d),$ $(d,d+l_2)$ (addition modulo~$n).$
\vskip 3pt
Proof.} Under the assumptions of the lemma, delays of at least one cluster occur on the limit cycle, and, at moment such that, at this moment a delay ends, the system is at one of the states $(l_1,0),$ $(0,l_2),$ $(d+l_1,d),$ $(d,d+l_2)$ (addition modulo~$n).$ This proves Lemma~2.
\vskip 3pt
{\bf Lemma 3.} {\it The system can be in a state of collapse if and only if
$$l_1>n-d.\eqno(3)$$
\vskip 3pt
Proof.} Taking into account the inequalities $l_1\le l_2$ and (3), we get the inequality $l_2>n- d.$ Therefore the states $(0,d)$ and $(d,0)$ are states of collapse.

If the system is in a state of collapse, then each cluster is at a node and occipies the other node. However, if (3) holds, a cluster cannot be at a node and occupy the other node simultaneously. This contradiction proves Lemma~3.

\section*{5. Criterion for resulting in a state of free movement}

\hskip 18pt  The following theorem and Lemma 1 provide a necessary and sufficient condition  of system resulting in free movement.
\vskip 3pt
{\bf Theorem 1.} {\it If the inequality
$$l_1+l_2\le n,\eqno(4)$$
holds, then the system results in a state of free movement from any initial state.

\vskip 3pt
Proof.} If (4) holds, then (3) does not hold, and, in accordance with Lemma~3, the system cannot be in a state of collapse. In accordance with
Lemma~2, it is sufficient to prove that from any of the states $(l_1,0),$ $(0,l_2),$ $(l_1+d,0),$ $(0,l_2+d)$ (addition modulo~$n),$ the system results in a state of free movement.

We shall consider all possible cases.
\vskip 3pt
1) Suppose $l_2\le d.$
$$A(t_0)=(l_1,0).$$
Then we have
$$A(t_0+d-l_1)=(d,d-l_1),\ A(t_0+d)=(d+l_1,d),$$
$$A(t_0+n-l_1)=(0,n-l_1),\ A(t_0+n)=(l_1,0).$$
Thus  the states $(l_1,0),$ $(d+l_1,d)$ are states of free movement.

We have the following sequences of states
$$A(t_0)=(0,l_2),\ A(t_0+d-l_2)=(d-l_2,d),$$
$$A(t_0+d)=(d,d+l_2),\ A(t_0+n-l_2)=(n-l_2,0),$$
$$A(t_0+n)=(0,l_2).$$
Hence the states $(0,l_2),$ $(d,d+l_2)$ are also states of free movement.
\vskip 3pt
2) Assume that $l_1\le d<l_2< n-d,$ $l_1+l_2\le n.$

We have
$$A(t_0)=(l_1,0),\ A(t_0+d-l_1)=(d,d-l_1),$$
$$A(t_0+d)=(d+l_1,d),\ A(t_0+n-l_1)=(0,n-l_1),$$
$$A(t_0+n)=(l_1,0);$$
$$A(t_0)=(0,l_2),\ A(t_0+d)=(d,d+l_2),$$
$$ A(t_0+n+l_2)=(n-l_2,0),\ A(t_0+n-l_2+d)=(n-l_2+d,d),$$
$$A(t_0+n)=(0,l_2).$$
Thus the states $(l_1,0),$ $(d+l_1,d),$  $(0,l_2),$ $(d,d+l_2)$
are states of free movement.
\vskip 3pt
3) If $l_1\le d <l_2,$ $l_2\ge n-d,$ $l_1+l_2\le n,$ then we have
$$A(t_0)=(l_1,0),\ A(t_0+d-l_1)=(d,d-l_1),$$
$$ A(t_0+d)=(d+l_1,d),\ A(t_0+n-l_1)=(0,1-l_1),$$
$$A(t_0+n)=(l_1,0);$$
$$A(t_0)=(0,l_2),\ A(t_0+n-l_2)=(n-l_2,0),$$
$$A(t_0+d)=(d,d+l_2-n),\ A(t_0+d+n-l_2)=(n-l_2+d,d),$$
$$A(t_0+n)=(0,l_2).$$
Whence, the states $(l_1,0),$ $(d+l_1,d),$  $(0,l_2),$ $(d,d+l_2-n)$
are states of free movement.
\vskip 3pt
4) Suppose $l_1> d,$ $l_1+l_2\le 1.$
 Then we have
$$A(t_0)=(l_1,0),\ A(t_0+d)=(l_1+d,d),$$
$$ A(t_0+n-l_1)=(0,n-l_1),\ A(t_0+n-l_1+d)=(d,1-l_1+d),$$
$$A(t_0+n)=(l_1,0),\ A(t_0)=(0,l_2),$$
$$A(t_0+d)=(d,d+l_2),\ A(t_0+n-l_2)=(n-l_2,0),$$
$$A(t_0+n)=(0,l_2).$$
Thus, the states $(l_1,0),$ $(d+l_1,d),$  $(0,l_2),$ $(d,d+l_2),$
are states of free movement.
\vskip 3pt
We have considered all possible cases under the assumption that the inequality $l_1+l_2\le n$ holds. Theorem~1 has been proved.

\section*{6. Behavior of the system under the conditions
$l_1\le n-d,$ $l_1+l_2>n$}

\hskip 18pt We shall prove the theorem that describes the behavior of the system in the case such that the system does not result in a state of free movement and collapse.
\vskip 3pt
{\bf Theorem 2.} {\it Assume that
$$l_1\le n-d,\ \  l_1+l_2>n.$$
Then, from any initial state, the system results in states of the same limit cycle with period
$$T=l_1+l_2,\eqno(5)$$
and the average velocity is equal to
$$v_1=v_2=\frac{n}{l_1+l_2}.\eqno(6)$$
\vskip 3pt
Proof.}  In accordance with Lemma~2, it is sufficient to prove that, from any initial states $(l_1,0),$ $(0,l_2),$ $(l_1+d,0),$ $(0,l_2+d)$ (addition modulo 1) the system results in a state of free movement.

We shall prove all possible cases.
\vskip 3pt
1) Suppose $l_1\le d,$ $l_1+l_2>1.$ Then we have the following sequences of transitions
$$A(t_0)=(l_1,0),\ A(t_0+d-l_1)=(d,d-l_1),$$
$$A(t_0+d)=(d+l_1,d),\ A(t_0+n-l_1)=(0,n-l_1),$$
$$A(t_0+l_2)=(0,l_2),\ A(t_0+n)=(n-l_2,0),$$
$$A(t_0+l_1+l_2)=(l_1,0),\ A(t_0)=(d,d+l_2-n),$$
$$A(t_0+n-l_2)=(d+n-l_2,d),\ A(t_0+l_1)=(d+l_1,d).$$
Thus the states $(l_1,0),$ $(0,l_2),$ $(d+l_1,d)$ belong to a limit circle with period $l_1+l_2$ and the average velocity of clusters calculated with the formulas (5), (6), and, from the state $(d,d+l_2-1),$ the system results in a state of this cycle.
\vskip 3pt
2) If $l_1> d,$ $l_2\le n-d,$ $l_1+l_2> n,$ then we have
$$A(t_0)=(l_1,0),\ A(t_0+d)=(d+l_1,d),$$
$$A(t_0+n-l_1)=(0, n-l_1),\ A(t_0+l_2)=(0,l_2),$$
$$ A(n)=(n-l_2,0),\ A(l_1+l_2)=(l_1,0),$$
$$A(t_0)=(0,d+l_2),\ A(t_0+n-d-l_2)=(n-l_2,0),$$
$$A(t_0+l_1-d)=(l_1,0).$$
Hence the states $(l_1,0),$ $(0,l_2),$ $(d+l_1,d)$ belong to a limit cycle with the period calculated with the formula (5) and the average velocity  calculated with the formula (6), and, from the state $(d,d+l_2-n),$ the system results in a state of this cycle.
\vskip 3pt
3) Assume that $d<l_1\le n- d,$ $l_2> n-d.$ In this case,
$$A(t_0)=(l_1,0),\ A(t_0+d)=(d+l_1,d),$$
$$ A(t_0+n-l_1)=(0,n-l_1),\ A(t_0+l_2)=(0,l_2),$$
$$A(t_0+n)=(n-l_2,0),\ A(t_0+l_1+l_2)=(l_1,0);$$
$$A(t_0)=(d,d+l_2-n),\ A(t_0+n-l_2)=(d+n-l_2,d),$$
$$ A(t_0+d)=(d+l_1,d).$$
Hence the states $(l_1,0),$ $(0,l_2),$ $(d+l_1,d)$ belong to a limit cycle with the period and average velocity calculated with the formulas (5), (6), and, from the state $(d,d+l_2-n),$ the system results in a state of this cycle.

Thus the states $(l_1,0),$ $(0,l_2),$ $(d+l_1,d)$ belong to a limit cycle with the period and average velocity of clusters calculated with the formulas (5), (6), and, from the state $(d,d+l_2-n),$ the system results in a state of this cycle.
\vskip 3pt
We have considered all possible cases under the assumption of Theorem~2. Thus Theorem~2 has been proved.

\section*{7. Criterion for resulting in a state of collapse}

The following theorem provides a necessary and sufficient condition  of system resulting in free movement

\vskip 3pt
{\bf Theorem  3.} {\it If (3) holds, then the system results in a state collapse from any initial state.
If (3) does not hold, then the system  does not result in a state collapse from any initial state.
\vskip 3pt
Proof.} If (3) holds, then (2) also holds, and, in accordance with Lemma~3, the system cannot be in a state of free movement. Therefore, in accordance with Lemma~2, it is sufficieent to prove that the system results in a state of collapse from any of the states
$(l_1,0),$ $(0,l_2),$ $(d+l_1-n,0),$ $(0,d+l_2-n).$

We have the sequences of states
$$A(t_0)=(l_1,0),\ A(n-l_1)=(0,n-l_1),\ A(d)=(0,d);$$
$$A(t_0)=(0,l_2),\ A(n-l_2)=(n-l_2,0),\ A(d)=(d,0);$$
$$A(t_0)=(d+l_1-n,0),\ A(n-l_1)=(d,0);$$
$$A(t_0)=(0,d+l_2-n),\ A(n-l_2)=(0,d).$$
Since the states $(0,d),$ $(d,0)$ are states of collapce, we get Theorem~3.

\begin{figure}[ht!]
\centerline{\includegraphics[width=250pt]{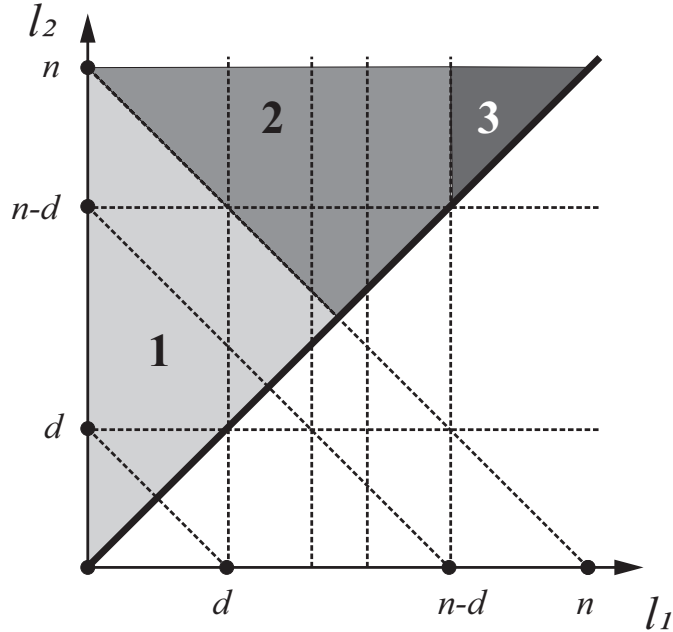}}
\caption{ 1) v=1, 2) $v =n/(l_1+l_2)$,  3) $v=0$}
\end{figure}

\section*{8. Conclusion}

\hskip 18pt We have proved that, for any $l_1,$ $l_2,$ $d,$ the average velocity of clusters is the same, i.e., $v_1=v_2,$ and the average velocity does not depend on the initial state. We have obtained the values of the average velocity. The values of the average velocity can be equal to 0, $n/(l_1+l_2)$ or 1.

The areas of velocity modes depending on $l_1$, $l_2$ are marked in fig. 2.

\end{document}